\documentclass[psfig,a4paper,11pt]{amsart}
\usepackage{amsfonts}
\theoremstyle{plain}
\begingroup
\newtheorem{theorem}{Theorem}[section]

\newtheorem{proposition}[theorem]{Proposition}

\newtheorem{remark}[theorem]{Remark}
\newtheorem{example}[theorem]{Example}
\newtheorem{definition}[theorem]{Definition}
\theoremstyle{definition}
\theoremstyle{remark}
\numberwithin{equation}{section}
\newcommand{\e}{\varepsilon}
\newcommand{\Om}{\Omega}

\newcommand{\dx}{\,dx}

\newcommand{\msim}{{\rm M}^{n\times n}_{\rm sym}} 
\newcommand{\lm}{{\mathcal L}^{n}}
\newcommand{\hn}{{\mathcal H}^{n-1}} 
\newcommand{\no}{\noindent}
 
\newcommand{\E}{{\mathcal E}}
\newcommand{\R}{{\mathbb R}}

\newcommand{\salt}{\noalign{\vskip .2truecm}}
 
\newcommand{\res}{\mathop{\hbox{\vrule height 7pt width .5pt depth 0pt
\vrule height .5pt width 6pt depth 0pt}}\nolimits}
\newcommand{\N}{{\mathbb N}}
 
\newcommand{\sbd}{{\rm SBD}(\Om)}
\newcommand{\bd}{{\rm BD}(\Om)}

\newcommand{\bom}{{\mathcal B}(\Om)}

\newcommand{\intcauchy}{\mskip 5mu -\mskip -18mu \int}
\newcommand{\B}{{\mathcal B}}

\newcommand{\norm}[1]{\left\Vert {#1} \right\Vert} 
\title[A lower semicontinuity result for some integral functionals in SBD]{A lower semicontinuity result for some \\
integral functionals in  
  the space SBD}
\author[Fran\c{c}ois Ebobisse]{Fran\c {c}ois Ebobisse}
\address{Department of Maths and Applied Maths, 
University of Cape Town, Rondebosch, 7700 South Africa}
\email{ebobisse@maths.uct.ac.za}

\begin{document}
\baselineskip3.2ex
\vskip .2truecm
\begin{abstract}
\small{The purpose of this paper is  to study the lower semicontinuity 
with respect to the strong $L^1$-convergence, of some integral functionals
 defined in the space ${\rm SBD}$ of 
special functions with bounded deformation. Precisely, we prove that, 
 if $u\in \sbd$,  $(u_h)\subset \sbd$ converges 
to $u$ strongly in $L^1(\Om,\R^n)$ and the mea\-sures $|E^ju_h|$ converge weakly $*$ to a measure 
$\nu$ singular with respect to the Lebesgue measure, then
$$\int_\Om f(x,\mathcal Eu)\dx\leq\liminf_{h\to\infty}\int_\Om f(x,\mathcal Eu_h)\dx$$
 provided $f$ satisfies some weak convexity property and the 
standard growth assumptions of order $p>1$.
}
\end{abstract}
\maketitle
{\small \tableofcontents}

\section{Introduction}
 Our goal in this paper is to extend in the framework of functions with bounded 
defor\-mation, the following lower semicontinuity theorem  
by Ambrosio \cite{AM1} for integral functionals defined in the space
SBV of special functions of bounded variation.
\begin{theorem}\label{theoamb}
Let $\Om\subset\R^n$ be an open set and let $f:\Om\times\R^k\times\R^{n\times k}$ 
be a Carath\'edory function satisfying:
\begin{itemize}
\item[(i)] for a.e. every $x\in\Om$, for every $(u,\xi)\in\R^k\times\R^{n\times k}$, 
$$ |\xi|^p\,\leq\,f(x,u,\xi)\,\leq\, a(x)+\Psi(|u|)(1+|\xi| ^p),$$ 
 where $p>1$, $a\in L^1(\Om)$ and the function $\Psi:[0,\infty)\to [0,\infty)$ is continuous;
\item[(ii)] for a.e. every $x\in\Om$ and every $u\in\R^k$, $f(x,u,\cdot)$ is quasi-convex.
\end{itemize}
 Then for every $u\in SBV(\Om,\R^k)$ and any sequence $(u_h)\subset SBV(\Om,\R^k)$ converging to $u$ 
 in $L^1_{loc}(\Om,\R^k)$ and such that 
\begin{equation}\label{jumpset} 
\sup_h\mathcal H^{n-1}(S_{u_h})<\infty
\end{equation}
 we have
$$\int_\Om f(x,u,\nabla u)\dx\leq\liminf_{h\to\infty}\int_\Om f(x,u_h,\nabla u_h)\dx.$$
\end{theorem}
Theorem \ref{theoamb} extends in the $SBV$ setting a classical lower semicontinuity 
result by Acerbi-Fusco \cite{AF} in the Sobolev space $W^{1,p}(\Om)$. 

 Later Kristensen in \cite{KRIS} extended Theorem \ref{theoamb} under the weaker assumptions  
\begin{equation}\label{kristjump}
\sup_h\int_{S_{u_h}}\theta(|u_h^+-u_h^-|)d\mathcal H^{n-1}<\infty
\end{equation}
for some function $\theta$ such that $\theta(r)/r\to\infty$ as $r\to 0^+$, 
 and $f$ is a normal integrand, i.e., for a.e. $x\in\Om$, $f(x,\cdot,\cdot)$ 
is lower semicontinuous in $\R^k\times\R^{n\times k}$ and there exists a Borel 
function $\tilde f:\Om\times\R^k\times\R^{n\times k}\to [0,\infty]$ 
such that $f(x,\cdot,\cdot)=\tilde f(x,\cdot,\cdot)$. \\
In the proof of 
Theorem \ref{theoamb} as well as in the Acerbi-Fusco result, the use of Lusin type 
approximation of functions in the given space ($BV$ or Sobolev spaces) 
by Lipschitz continuous functions is crucial.

Recently, Theorem \ref{theoamb} has been extended by 
Fonseca-Leoni-Paroni \cite{FOLEPA} to functionals depending also on the hessian matrices.

In this paper we deal with first order variational problem, 
but with integral functionals depending explicitely on the symmetrized 
derivative $Eu:=(Du+Du^T)/2$ 
and defined in the space $BD$ of functions with bounded deformation.

The main result of the paper is  the following lower semicontinuity 
theorem:
\begin{theorem}\label{lower1}
Let $p>1$ and Let $f:\Om\times\msim\to [0,\infty)$ be a Carath\'eodory function satisfying:
\begin{itemize}
\item[(i)] for a.e. every $x\in\Om$, for every $\xi\in \msim$, 
$$\displaystyle\frac{1}{C}|\xi|^p\,\leq\,f(x,\xi)\,\leq\, \phi(x)+C(1+|\xi| ^p),$$ 
for some constant $C>0$ and a function $\phi\in L^1(\Om)$;
\item[(ii)] for a.e. every $x_0\in\Om$, $f(x_0,\cdot)$ is symmetric 
quasi-convex i.e., 
\begin{equation}\label{symqc}
f(x_0,\xi)\leq\intcauchy _A f(x_0,\xi+{\E}\varphi (x))dx
\end{equation}
for every bounded open subset $A$ of $\R ^n$, for every $\varphi \in W^{1,\infty }_0(A,\R ^n )$ and $\xi\in\msim$.
\end{itemize}
Then for every $u\in \sbd$, for any sequence $(u_h)\subset \sbd$ 
converging to $u$ strongly in $L^1(\Om,\R^n)$
with $|E^ju_h|$ converging weakly $*$ to a positive 
measure $\nu$ singular with respect to the Lebesgue measure, 
we have
\begin{equation}\label{lower2}
\nonumber\int_\Om f(x,\E u)\,dx\,\leq\,\liminf_{h\to\infty}\int_\Om f(x,\E u_h)\dx.
\end{equation}
\end{theorem}
In the literature there are various results  
on lower semicontinuity and relaxation
 of convex integral functionals in BD with linear growth in the strain tensor, in 
connection with elasto-plasticity problems (see \cite{AG, SUQ, TEM, BDV}). 
Concerning non convex functionals with linear growth
 we mention the papers \cite{EBO2, BFT, ET}. As far as the author
 knows, there is no result on lower 
semicontinuity of non convex volume energies with 
superlinear growth in the strain tensor. So, Theorem \ref{lower1} is the first
 lower semicontinuity result for this class of functionals.

The proof of Theorem \ref{lower1} follows the lines of Theorem \ref{theoamb}.
 We use the blow-up method introduced in \cite{FOMU1} and described as a 
two-steps process whose first step 
here is the proof of a lower semicontinuity result whenever $\Om$ is the unit
ball $B(0,1)$, the limit function
 is linear and  $|E^ju_h|(B(0,1))$ converge to zero 
(see Proposition \ref{step1}). 
In a second step, we use a blow-up argument
 through the approximate 
differentiability of BD functions to reduce the problem into the first 
step.

The use of Lusin type approximation of BD functions by Lipschitz 
functions is crucial 
in the proof of  Proposition \ref{step1}. This result  established in 
\cite{EBO1} and refined here in Proposition \ref{lusbdp} is obtained using
 a "Poincar\'e type" inequality for BD functions (see Theorems 
\ref{kohn1} and \ref{kohn2}) together with the maximal 
function of Radon measures.

This paper is organized as follows. In section 2 we collect and prove some fine
 properties of BD functions that will be used in the proof of our main 
result. Section 3 is devoted to the proof of Theorem \ref{lower1}. 
  In section 4, we discuss the assumption (in Theorem \ref{lower1}) that
 the measures $|E^ju_h|$ converge weakly $*$ to a positive measure $\nu$
 singular with respect to the Lebesgue measure. In Example \ref{exp0}, we 
consider a minimization problem in SBD with a unilateral constraint on the 
jump sets and we show that  minimizing sequences satisfy the assumption
 on $|E^ju_h|$. However, this 
assumption is not always compatible with the SBD compactness criterion 
 (Theorem \ref{Compactness-SBD}). In fact, we construct in Example 
\ref{examp1}, a sequence of functions $(u_h)$ which verifies 
the assumptions of 
Theorem \ref{Compactness-SBD} while $|E^ju_h|$ converge 
weakly $*$ to a measure proportional to the Lebesgue measure. 
 
\section{Notation and preliminaries}

Let $n\geq 1$ be an integer. 
We denote by $M^{n\times n}$ the space of $n\times n$ matrices and by 
$\msim $ the subspace of symmetric matrices in $M^{n\times n}$. For any 
$\xi\in M^{n\times n}$, $\xi ^T$ is the transpose of $\xi$. 
Given $u,\,v\in\R ^n$, $u\otimes v$ and $u\odot v:=(u\otimes v+v\otimes u)/2$ 
denote the tensor and symmetric products of $u$ and $v$, respectively. 
We use the standard notation, $\lm $ and $\hn $ to denote respectively the 
Lebesgue outer measure and the $(n-1)$-dimensional Hausdorff measure. 
For every
set $E\subset\R^n$,\, $\overline E$, $\,|E|$ and $\chi_E$ stand respectively for the closure 
of $E$, the Lebesgue outer measure of $E$ and the charateristic function of $E$, that is
$\chi_E(x)=1$ if $x\in E$ and  $\chi_E(x)=0$ if $x\notin E$.
For $1\leq p\leq\infty$, 
$\norm{\cdot}_p$ will denote the norm in the $L^p$ space.

Let $\Om $ be an open subset of $\R ^n$ We denote by  $\bom$ 
the family of Borel subsets of $\Om$.
For any $x\in \Om $ and $\rho>0$,  $B(x,\rho )$ denotes the open ball of 
$\R^n $ centered at
$x$ with radius $\rho $.
When $x=0$ and $\rho=1$ we simply write $B_1$. We use the notation $w_n$ for
 the Lebesgue measure of the ball $B_1$.
If $\mu $ is a Radon measure, we denote  $|\mu |$ its total variation.

\begin{definition}\label{bd}
A function $u:\Om \to \R^n $ is with {\it bounded deformation} in $\Om $ 
if $u\in L^1(\Om ,\R^n)$ and 
$Eu:=(Du+Du^T)/2\in \mathcal M_b\bigl(\Om ,\msim\bigr)$, 
where $Du$ is the distributional gradient of $u$ and 
$\mathcal M_b\bigl(\Om,\msim\bigl)$ 
is the space of $\msim$-valued
Radon measures with finite total variation in $\Om $.
\end{definition}
The space $\bd$ of functions with bounded deformation in $\Om$ was introduced
 in \cite{MSC} and
studied, for instance in \cite{AG}, \cite{KH},
 \cite{SUQ}, \cite{TEM} in relation with the static model of Hencky 
in perfect plasticity.
$\bd$ is a Banach space when equipped with the norm
$$
\norm{u}_{BD(\Om )}:=\norm{u}_{L^1(\Om ,\R ^n)}+\vert Eu\vert (\Om)
$$
where $|Eu|(\Om)$ is the total variation of the measure $Eu$ in $\Om$. \\
Whenever the open set $\Om $ is assumed to be connected, the kernel of the 
operator $E$ is the class of {\it rigid motions} denoted here by 
${\mathcal R}$, and composed of affine maps of the form $Mx+b$, where $M$ is a 
skew-symmetric $n\times n$ matrix and $b\in\R^n$. Therefore ${\mathcal R}$ is 
closed and finite-dimensional.\\
Fine properties of ${\rm BD}$ functions were studied, for instance, in 
\cite{ACDM}, \cite{BCDM} and \cite{KH}. The following ``Poincar\'e type''
inequality for $BD$ functions has been proved by Kohn \cite{KH} (see also \cite{ACDM}).
\begin{theorem}\label{kohn1}
{\it Let $\Om $ be a bounded connected open subset of $\R ^n$ with
  Lipschitz boundary. Let $R:BD(\Om )\to {\mathcal R}$ be a continuous
  linear map which leaves ${\mathcal R}$ fixed.\\
Then there exists a positive constant $C(\Om ,R)$ such that:
\begin{equation}\label{poincare'}
 \int _\Om \vert u-R(u)\vert dx\leq C(\Om ,R)\vert Eu\vert (\Om
)\quad \mbox { for any }u\in BD(\Om ).
\end{equation} }
\end{theorem}
When $\Om $ is an open ball of $\R ^n $ the\-re is a precise representation
of the rigid motion $R(u)$, given in the following theorem.
\begin{theorem}\label{kohn2}
{\it Let $u\in BD(\R ^n )$, $x\in \R ^n $ and $\rho >0$. 
Then there exists a vector $d_\rho (u)(x)\in \R ^n$ and 
an $n\times n$ skew-symmetric matrix $A_\rho (u)(x)$ such that:
\begin{equation}\label{precisepoincare'}
\int _{B(x,\rho)}\vert u(y)-d_\rho (u)(x)-A_\rho (u)(x)(y-x)\vert
dy\leq C(n)\rho \vert Eu\vert (B(x,\rho))
\end{equation}
where $C(n)$ is a positive constant depending only on the dimension
$n$.\\
Moreover, $d_\rho (u)(x)$ and $A_\rho (u)(x)$ are expressed as singular
integrals in the following ways:
\begin{equation}\label{vect}
 d^i_\rho (u)(x):=\sum \limits _{l,m=1}^n\int _{\vert y-x\vert \geq \rho
  }{{\bigwedge ^i_{lm}(y-x)}\over {nw_n\vert y-x\vert
  ^n}}dEu_{lm}(y);
\end{equation}
\begin{equation}\label{anti}
A^{ij}_\rho (u)(x):=\sum \limits _{l,m=1}^n\int _{\vert y-x\vert \geq
  \rho }-{{\Gamma ^{lm}_{ij}(y-x)}\over {2w_n\vert y-x\vert
  ^{n+2}}}dEu_{lm}(y),
\end{equation}
where $\bigwedge $ and $\Gamma $,  respectively third and fourth-order 
tensor valued functions, are defined and studied in \cite{KH}, \cite{ACDM}.}
\end{theorem}

We recall that if $u\in\bd$, then the 
jump set $J_u$ of $u$ is a countably 
$(\hn,n-1)$-rectifiable Borel set and the 
following decomposition of the measure $Eu$ holds
\begin{equation}\label{decomp}
Eu=\E u\lm\, +\, E^su=\E u\lm\, +\,E^ju\,+\,E^cu\,,
\end{equation}
where $E^ju:=([u]\odot\nu_u)\hn\res J_u$, $[u]:= u^+-u^-$, 
$u^+$ and $u^-$ are the {\it one-sided Lebesgue limits} of $u$ with 
respect to the measure theoretic normal $\nu_u$ of $J_u$, 
$\E u$ is the density of the absolutely continuous part of $Eu$ with respect 
to $\lm$, $E^su$ is the singular part, and $E^cu$ is the {\it Cantor part} 
and vanishes on the Borel sets that are $\sigma$-finite with respect to $\hn$ 
(see \cite{ACDM}). 
\vskip .2truecm
 Hereinafter we will use the following Proposition proved
in  \cite[Proposition 7.8 and Remark 7.9]{ACDM}.
\begin{proposition}\label{singinte}
Let $K:\R ^n\setminus\{ 0\}\to \R $ be a $0$-homogeneous
  function, smooth and with mean value zero on the unit sphere
  ${\bf S}^{n-1}$. For any Radon measure $\mu $ with finite total variation
  in $\R ^n $, let us define the functions 
$$h_\rho (x):=\int _{\vert y-x\vert \geq \rho }{{K(y-x)}\over {\vert
    y-x\vert ^n}}d\mu (y)\quad \quad \rho >0.$$
Then the function $h(x):=\sup \limits _{\rho >0}\vert h_\rho (x)\vert $
    satisfies the following weak $L^1$ estimate
\begin{equation}\label{weakestimate1}
\bigl|\{x\in \R ^n \mbox {: }h(x)>t\}\bigr|\leq {{C(n,K)}\over
  t}\vert \mu \vert (\R ^n ).
\end{equation} 
Moreover, if $\mu=f\mathcal L^n$ with $f\in L^p(\R^n)$, 
then the following strong $L^p$ estimate holds
\begin{equation}\label{strongestimate1}
\norm{h}_p\leq C(n,K)\norm{f}_p.
\end{equation}
\end{proposition}

Let us also recall the theorem by Ambrosio-Coscia-Dal Maso \cite{ACDM}
on the approximate differentiability of ${\rm BD}$ functions.
\begin{theorem}\label{approx}
Let $\Om $ be a bounded open set in $\R ^n$ with Lipschitz boundary. Let 
$u\in\bd$. Then for $\lm$ almost every 
$x\in\Om $ there exists an $n\times n$ matrix 
$\nabla u(x)$ such that
\begin{equation}\label{apdif}
\lim_{\rho\to0}\frac{1}{\rho ^n}\int_{B(x,\rho)}\frac{|u(y)-u(x)-
\nabla u(x)(y-x)|}{\rho}dy=0\,,
\end{equation}
and 
\begin{equation}\label{symapp}
\lim_{\rho\to0}\frac{1}{\rho ^n}\int_{B(x,\rho)}
\frac{|(u(y)-u(x)-\E u(x)(y-x),y-x)|}{|y-x|^2}dy=0
\end{equation}
for $\lm$-almost every $x\in\Om $.
\end{theorem} 
\no In particular, from (\ref{apdif}) we have $u$ is approximately differentiable 
$\lm$-almost everywhere in 
$\Om$ and from Proposition \ref{singinte} it has been proved the function $\nabla u$
 satisfies the weak $L^1$ estimate
$$
\lm (\{x\in\Om\mbox{ : } |\nabla u(x)|>t\})\leq
\frac{C(n,\Om)}{t}\norm{u}_{\bd}\quad\forall t>0,
$$
where $C(n,\Om)$ is a positive constant depending only on $n$ and $\Om $.

\no From (\ref{symapp}) and (\ref{apdif}) one can easily see that 
\begin{equation}\label{symappdiff}
\E u(x)=(\nabla u(x)+\nabla u(x)^T)/2\quad\mbox{for }\lm \mbox{-a.e. }x\in\Om.
\end{equation}

\no Analogously to the space ${\rm SBV}$ introduced by De Giorgi and 
Ambrosio (see for instance \cite{AMFUPA}), the space ${\rm SBD}$ was introduced by Bellettini and Coscia in \cite{BC} and studied  in \cite{BCDM}.
\begin{definition}\label{sbd} 
The space $\sbd$ of {\em special functions with bounded deformation}, is the 
space of functions $u\in\bd$ such that the measure $E^cu$ in (\ref{decomp}) 
is zero. 
\end{definition}
We set 
$$A(u)(x):=\sup_{\rho>0}|A_\rho (u)(x)|$$
with $A_\rho (u)$ the anti-symmetric matrix defined in (\ref{anti}). Note that for every
$u\in SBD(\Om)$,  $A_\rho (u)=L_\rho(u)+J_\rho(u)$ with
\begin{equation}\label{leb}
L^{ij}_\rho (u)(x):=\sum \limits _{l,m=1}^n\int _{\vert y-x\vert \geq
  \rho }-{{\Gamma ^{lm}_{ij}(y-x)}\over {2w_n\vert y-x\vert
  ^{n+2}}}{\mathcal E}u_{lm}(y)dy
\end{equation}
and 
\begin{equation}\label{jump}
J^{ij}_\rho (u)(x):=\sum \limits _{l,m=1}^n\int _{\vert y-x\vert \geq
  \rho }-{{\Gamma ^{lm}_{ij}(y-x)}\over {2w_n\vert y-x\vert
  ^{n+2}}}dE^ju_{lm}(y).
\end{equation}
We set also
\begin{equation}\label{abscont}
L(u)(x):=\sup_{\rho>0}|L_\rho (u)(x)|\quad\mbox{ and }
\quad J(u)(x):=\sup_{\rho>0}|J_\rho (u)(x)|.\end{equation}
Let us also recall that, given a $\R^m$-valued Radon Measure $\mu$ in $\R^n$, 
the maximal function of $\mu$ is defined by
$$M(\mu )(x):=\sup \limits _{\rho >0}{\vert \mu \vert 
(B(x,\rho))\over \left|B(x,\rho)\right|}\quad\quad \forall x\in \R ^n
    .$$
Whenever $\mu=g\mathcal L^n$, we recover the maximal function of the function
$g$ (see \cite{STN}).
\vskip .3truecm
 The following theorem is proved in \cite{EBO1}.
\begin{theorem}\label{lusbd}
Let $\Om $ be either $\R ^n$ or a Lipschitz bounded 
open subset of $\R ^n $ and $u\in BD(\Om )$.
Then for any $\lambda >0$, there exists a  Lipschitz continuous function
$v_\lambda :\Om \to \R ^n $
with lip$(v_\lambda ) \leq C\lambda $ such that:
\begin{equation}\label{exept1}
\bigl|\{x\in \Om \mbox { : }v_\lambda (x)\neq u(x)\}\bigr|\leq
  {C\over
  \lambda }\norm{u} _{BD(\Om )}, 
\end{equation}
where $C$ is a positive constant only depending on $n$ or also on
$\Om $.

\end{theorem}
\no We can further refine the estimate (\ref{exept1}) when the func\-tion 
$u\in SBD(\Om)$ with $\mathcal Eu\in L^p(\Om,\msim)$.
\begin{proposition}\label{lusbdp}
Let $p\in (1,\infty)$ and $u\in SBD(\R^n)$ with 
$\mathcal Eu\in L^p(\R^n,\msim)$. For every $\lambda>0$ 
 there exists a function $v_\lambda :\R^n\to\R^n$ 
Lipschitz continuous with lip$(v_\lambda)\leq C\lambda$, 
$|v_\lambda(x)|\leq C\lambda\quad\mbox{ 
for every }x\in\R^n$,  
and for any Borel subset $E$ of $\R^n$, the following estimate holds
\begin{eqnarray}\label{exep2}
 & { } &\bigl|E\cap\{x\in \R^n \mbox { : }v_\lambda (x)\neq u(x)\}\bigr|  \leq 
  {C\over
  \lambda }\bigl[||u||_{L^1(\R^n,\R^n)}+|E^ju|(E)\bigr] + \\
\salt
\nonumber & { } & \hskip 2truecm +\,{1\over
  \lambda ^p}\int_{E\cap\{L(u)(x)>\lambda\}}|L(u)(x)|^pdx +{1\over
  \lambda ^p}\int_{E\cap\{M(|\E u|)>\lambda\}}[M(|\E u|)]^pdx.
\end{eqnarray}
where $C$ is a positive constant only depending on $n$.
\end{proposition}
\begin{proof}
For $\lambda>0$, we set 
$$E_\lambda:=\{x\in\R^n\mbox{: } M(|u|\lm +|Eu|)(x)\leq 3\lambda\mbox{ and }A(u)(x)\leq 2\lambda\}.$$
It has been proved in Theorem \ref{lusbd} that $u|_{E_\lambda\setminus S_u}$ 
is Lipschitz continuous with Lipschitz constant less or equal to a positive 
constant proportional to $\lambda$.  Moreover, from Lebesgue 
differentiation theorem we have also
$$|u(x)|\leq 3\lambda\quad\forall x\in E_\lambda\setminus S_u.$$
The function $v_\lambda$ is then obtained 
from $u|_{E_\lambda\setminus S_u}$ by 
Kirszbraun's Theorem (see Federer \cite[Theorem 2.10.43]{FD}). \\
Now given $E\in {\mathcal B}(\R^n)$, since 
$E\cap\{x\in \R^n \mbox { : }v_\lambda (x)\neq u(x)\}\,
\subset\,E\setminus (E_\lambda\setminus S_u)$,  
it is sufficient to estimate the measure of the Borel set 
$E\setminus (E_\lambda\setminus S_u)$.\\
Note that 
\begin{eqnarray*}
& &\bigl|E\cap\{x\in\R^n\mbox{: }A(u)(x)>2\lambda\}\bigr|
\leq\bigl|E\cap\{x\in\R^n\mbox{: }L(u)(x)>\lambda \}\bigr|\\
\salt
& &\hskip 5.4truecm  +\,\,
\bigl|E\cap\{x\in\R^n\mbox{: }J(u)(x)>\lambda \}\bigr|
\end{eqnarray*}
where $L$ and $J$ are defined in (\ref{abscont}). 
 From Proposition \ref{singinte} and Chebychev inequality we get
 respectively
$$\bigl|E\cap\{x\in\R^n\mbox{: }J(u)(x)>\lambda\}\bigr|\leq {C(n)\over \lambda}|E^ju|(E)$$
and
$$\bigl|E\cap\{x\in\R^n\mbox{: }L(u)(x)>\lambda \}\bigr|
\leq {1\over \lambda^p}\int_{E\cap\{L(u)(x)>\lambda\}}|L(u)|^pdx.$$
So, we obtain
\begin{equation}\label{anti1}
\bigl|E\cap\{x\in\R^n\mbox{: }A(u)(x)>2\lambda\}\bigr|
\leq {C(n)\over \lambda}|E^ju|(E)\,+\,{1\over
  \lambda ^p}\int_{E\cap\{L(u)(x)>\lambda\}}|L(u)|^pdx.
\end{equation}
On the other hand, using covering theorems (see \cite{AMFUPA}, \cite{FD}) 
and the properties of maximal functions of $L^p$ functions, we obtain 
the estimates
\begin{eqnarray}\label{maxi1}
& &\bigl|E\cap\{x\in\R^n\mbox{: }M(|u|\lm +|Eu|)(x)>3\lambda\}\bigr|\,\leq \,
\bigl|E\cap\{x\in\R^n\mbox{: }M(|u|\lm)(x)>\lambda\}\bigr|\\
\salt
\nonumber & &\hskip 2.truecm\,+\,\bigl|E\cap\{x\in\R^n\mbox{: }M(E^ju)(x)>\lambda\}\bigr|
\,+\, 
\bigl|E\cap\{x\in\R^n\mbox{: }M(\E u)(x)>\lambda\}\bigr|\\
\salt
\nonumber & &\hskip 1truecm\,\leq \, {C(n)\over \lambda}\bigl[||u||_{L^1(\R^n,\R^n)}+
|E^ju|(E)\bigr]\,+\,{ 1\over
  \lambda ^p}\int_{E\cap\{M(\E u)>\lambda\}}[M(\E u)]^p\dx.
\end{eqnarray}
Therefore the estimate (\ref{exep2}) is obtained by 
 adding (\ref{anti1}) to (\ref{maxi1}).
\end{proof}
\begin{remark}\label{lusindom}
{\rm Let $\Om$ be a bounded connected open subset of $\R^n$ 
with Lipschitz continuous boundary $\partial\Om$. Let $u\in SBD(\Om)$ with 
$\mathcal Eu\in L^p(\Om,\msim)$. Let $\overline u$ be  the extension
 of $u$ by $0$ outside $\Om$. We recall that
$$E\overline u=\mathcal Eu\mathcal L^n\res\Om\,+\,E^ju\res\Om\,-\,
\gamma(u)\odot\nu\mathcal H^{n-1}\res\partial\Om$$
where $\gamma(u)$ and $\nu$ are respectively the trace of $u$ on 
$\partial\Om$ and the outer unit normal vector  to $\partial\Om$. 
So, applying Proposition \ref{lusbdp}  to $\overline{u}$, we get the 
following estimate 
\begin{eqnarray}\label{exep13}
& &\bigl|E\cap\{x\in \R^n \mbox { : }v_\lambda (x)\neq\overline{u}(x)\}\bigr|\\
\salt
\nonumber & \leq & 
  {C(n,\Om)\over
  \lambda }\left[||u||_{L^1(\Om,\R^n)}+|E^ju|(\Om\cap E)+
\int_{\partial\Om\cap E}|\gamma(u)\odot\nu|d{\mathcal H}^{n-1}
\right]\\
\salt
\nonumber & & +\, {1\over
  \lambda ^p}\int_{E\cap\{L(u)\,>\,\lambda\}}|L(u)|^pdx
\, +\,{1\over
  \lambda ^p}\int_{E\cap\{M(|\E u|)>\lambda\}}[M(|\E u|)]^p\dx
\end{eqnarray}
 for every $E\in\mathcal B(\R^n)$. 
In particular 
 for any $E\in\mathcal B(\Om)$ we get
\begin{eqnarray}\label{exep3}
& &\bigl|E\cap\{x\in \Om \mbox { : }v_\lambda (x)\neq u(x)\}\bigr|
\,\,\leq \,\, 
  {C(n,\Om)\over
  \lambda }\left[||u||_{L^1(\Om,\R^n)}+|E^ju|(\Om)\right]\\
\salt
\nonumber & &\hskip 1truecm\,+\,{1\over
  \lambda ^p}\int_{E\cap\{L(u)(x)\,>\,\lambda\}}|L(u)|^pdx
\, +\,{1\over
  \lambda ^p}\int_{E\cap\{M(|\E u|)>\lambda\}}[M(|\E u|)]^p\dx.
\end{eqnarray}

}
\end{remark}
\section{The proof of the main result}
 This section is essentially devoted to the proof of Theorem  \ref{lower1} 
where the following propo\-si\-tion will be crucial.
\begin{proposition}\label{step1}
Let $f_h:\Om\times\msim\to [0,\infty)$ be a 
sequence of Carath\'eodory functions satisfying 
for a.e. every $x\in\Om$, for every $\xi\in \msim$, 
$$\displaystyle\frac{1}{C}|\xi|^p\,\leq\,f_h(x,\xi)
\,\leq\, \phi_h(x)+C(1+|\xi| ^p),$$ 
for some constant $C>0$ and a sequence $(\phi_h)$ uniformly bounded 
in $L^1(B_1)$.
Assume that there exist an $\lm$-negligible 
set $N\subset B_1$ and a symmetric quasi-convex function 
$f:\msim\to [0,\infty)$ such that 
$\lim_{h\to\infty}f_h(y,\xi)=f(\xi)$ uniformly on compact subsets of $\msim$ 
and for any $y\in B_1\setminus N$.
Then, for any sequence $(u_h)$ in SBD$(B_1)$ converging strongly 
in $L^1(B_1,\R^n)$ to a linear function $u$, with $\lim_{h\to\infty}|E^ju_h|(B_1)\to 0$, we have
$$
\int_{B_1} f(\E u)\,dx\,\leq\,\liminf_{h\to\infty}\int_{B_1} f_h(x,\E u_h)\dx.
$$
\end{proposition}
\begin{proof}
Let $(u_h)\subset SBD(B_1)$ be a sequence which converges strongly in 
$L^1(B_1,\R^n)$ to a linear function $u$ and 
$\lim_{h\to\infty}|E^ju_h|(B_1)\to 0$. Up to substituting $u_h$ by $u_h-u$ and $f_h(x,z)$
by $f_h(x, z+\E u)$ we can assume that
$u\equiv 0$. So, we have to prove that 
\begin{equation}\label{semiinf0}
|B_1| f( 0)\,\leq\,\liminf_{h\to\infty}\int_{B_1} f_h(x,\E u_h)\dx.
\end{equation}
Up to a subsequence we assume that 
$$\liminf_{h\to\infty}\int_{B_1} f_h(x,\E u_h)\dx
=\lim_{h\to\infty}\int_{B_1} f_h(x,\E u_h)\dx<\infty.$$
So the sequence $(\E u_h)$ is  uniformly bounded in $L^p(B_1,\msim)$.
We set 
\begin{equation}\label{psi}
\Psi_h:=[M(\E u_h)]^p+|L(u_h)|^p +|\phi_h|
\end{equation}
where $M$ is the maximal function and $L$ is defined in (\ref{abscont}).
From the assumptions and from Proposition \ref{singinte}, we have that $(\Psi_h)$ is a 
bounded sequence in $L^1(B_1)$. So, By Chacon Bitting Lemma 
(see for instance \cite[Lemma 5.32]{AMFUPA}) there 
exist a subsequence of $(\Psi_h)$ (still denoted $(\Psi_h)$)
and a decreasing sequence of sets $(E_k)\subset\mathcal B(B_1)$ 
such that $|E_k|\to 0$ as $k\to\infty$ and the sequence
$(\Psi_h1_{B_1\setminus E_k})_h$ is equiintegrable for any $k\in\N$. 
We introduce the following modulus of equiintegrability
for the sequence $(\Psi_h1_{B_1\setminus E_k})_h$
\begin{equation}\label{modeq}
W_k(\delta):=\sup\left\{\begin{array}{ll}
\displaystyle\limsup_{h\to\infty}\int_F\Psi_h\dx\mbox{: }  
F\in\mathcal B(B_1)\mbox{ with} \\
\salt
F\subset B_1\setminus E_k\mbox{ and }|F|\leq\delta 
\end{array}
\right\}\quad\forall\delta >0,\,\forall k\in\N.
\end{equation}
It follows that $W_k(\delta)\to 0$ as $\delta\to 0$.\\
Now, from Proposition \ref{lusbdp} and Remark \ref{lusindom} we have for 
any integer $m\geq 1$, 
a Lipschitz continuous function $v_{h,m}: \overline B_1\to\R^n$
 and a set $E_{h,m}\in\mathcal B(B_1)$  such that 
\begin{equation}\label{lipco}
Lip(v_{h,m})\leq C(n)m,\quad |v_{h,m}(x)|\leq C(n)m\quad
\forall x\in \overline {B_1},\quad v_{h,m}=u_h\mbox{ in }B_1\setminus E_{h,m}
\end{equation}
 and
for any Borel subset $E$ of $B_1$ the following estimate holds
\begin{eqnarray}\label{estim1}
& &\bigl|E_{h,m}\setminus E\bigr| \leq 
  {C(n)\over m }\bigl[||u_h||_{L^1(B_1,\R^n)}+|E^ju_h|(B_1)\bigr]\\
\salt
\nonumber & &
\,+\,{1\over m^p}\int_{\{L(u_h)(x)>m\}\setminus E}|L(u_h)(x)|^pdx
\,+\,{1\over
  m ^p}\int_{\{M(|\E u_h|)>m\}\setminus E}[M(|\E u_h|)]^pdx.
\end{eqnarray}
In particular for $E=E_k$ we get from the definition of $\Psi_h$ that
\begin{equation}\label{estim2}
\bigl|E_{h,m}\setminus E_k\bigr| \leq {C(n)\over m }
\bigl[||u_h||_{L^1(B_1,\R^n)}+|E^ju_h|(B_1)\bigr] +
{2\over m^p}\int_{\{\Psi_h>m^p\}\setminus E_k}\Psi_h\dx.
\end{equation}
We set $S:=\sup_h||\Psi_h||_1$. Using the fact that 
$\bigl|\{\Psi_h>m^p\}\bigr|\leq {S\over m^p}$ together with $u_h\to 0$ strongly in 
$L^1(B_1,\R^n)$ and $|E^ju_h|(B_1)\to 0$ (by assumptions), 
we get from (\ref{estim2})
that
\begin{equation}\label{estim3}
\limsup_{h\to\infty}\,m^p\bigl|E_{h,m}\setminus E_k\bigr| 
\leq 2W_k\Bigl({S\over m^p}\Bigr).
\end{equation}
From the inequality (\ref{estim3}), it is easy to see (for $m$ large enough) 
that 
$$\limsup_{h\to\infty}\int_{E_{h,m}\setminus E_k}\phi_h\dx\leq
\limsup_{h\to\infty}\int_{E_{h,m}\setminus E_k}\Psi_h\dx\leq 
2W_k\Bigl({S\over m^p}\Bigr).$$
Now from (\ref{lipco}), it follows by Ascoli-Arzel\`a that the sequence
 $(v_{h,m})_h$ is relatively compact in $C(\overline B_1,\R^n)$. 
Hence, we get up to a subsequence, that for every integer $m\geq 1$, 
$v_{h,m}$ converges uniformly to
a function $v_m\in C(\overline B_1,\R^n)$ as $h\to\infty$. 
Since, $|E_k|\to 0$ as $k\to\infty$, to get (\ref{semiinf0}), it is enough to prove that
\begin{equation}\label{seminf1}
|B_1\setminus E_k| f( 0)\,\leq\,\liminf_{h\to\infty}\int_{B_1} f_h(x,\E u_h)\dx
\qquad\forall k\in\N.
\end{equation}
We have the following estimates
\begin{eqnarray*}
& &\int_{B_1} f_h(x,\E u_h)\dx
\geq\int_{B_1\setminus (E_{h,m}\cup E_k)} f_h(x,\E u_h)\dx
\,=\,\int_{B_1\setminus (E_{h,m}\cup E_k)} f_h(x,\E v_{h,m})\dx\\
\salt
& & \hskip 3truecm =\,\int_{B_1\setminus  E_k} f_h(x,\E v_{h,m})\dx
\,-\,\int_{E_{h,m}\setminus E_k} f_h(x,\E v_{h,m})\dx\\
\salt
& & \hskip 3truecm \geq\,\int_{B_1\setminus  E_k} f_h(x,\E v_{h,m})\dx
\,-\,\int_{E_{h,m}\setminus E_k}\phi_hdx -Cm^p\bigl|E_{h,m}\setminus E_k\bigr|.
\end{eqnarray*}
So, passing to the limit as $h\to\infty$, and using (\ref{estim2}) 
and  (\ref{estim3}) we get that
\begin{equation}\label{seminf2}
\liminf_{h\to\infty}\int_{B_1} f_h(x,\E u_h)\dx\geq
\liminf_{h\to\infty}\int_{B_1\setminus  E_k} f_h(x,\E v_{h,m})\dx
\,-\,CW_k\Bigl({S\over m^p}\Bigr).
\end{equation}
From the assumption of the convergence of $f_h(x,\xi)$ to $f(\xi)$, we get that 
\begin{equation}\label{seminf3}
\liminf_{h\to\infty}\int_{B_1\setminus  E_k} f_h(x,\E v_{h,m})\dx
\geq\liminf_{h\to\infty}\int_{B_1\setminus  E_k} f(\E v_{h,m})\dx.
\end{equation}
Using the symmetric quasi-convexity of the function $f$, we get also that
\begin{equation}\label{seminf4}
\liminf_{h\to\infty}\int_{B_1\setminus  E_k} f(\E v_{h,m})\dx\geq \int_{B_1\setminus  E_k} f(\E v_m)\dx.
\end{equation}
Indeed, $f$ symmetric quasi-convex means that $f\circ\pi$ is quasi-convex
 in the classical sense, where $\pi$ is the projection on symmetric matrix.
 Since $lip(v_{h,m})\leq C(\Om, n)m$, it is easy to see that the $(v_{h,m})_h$
 converges  weakly $\star$ in $W^{1,\infty}(B_1,\R^n)$ to the function $v_m$
 and hence  (\ref{seminf4}) follows from a classical lower semicontinuity theorem 
by Morrey (see for instance Dacorogna \cite{dac}).\\
Finally from (\ref{seminf2}), (\ref{seminf3}) and  (\ref{seminf4}) we obtain that
\begin{equation}\label{seminf5}
\liminf_{h\to\infty}\int_{B_1} f_h(x,\E u_h)\dx
\geq \int_{B_1\setminus  E_k} f(\E v_m)\dx-\,CW_k\Bigl({S\over m^p}\Bigr).
\end{equation}
On the other hand, from (\ref{estim3}) we have also that
\begin{equation}\label{estim4}
m^p\bigl|\{x\in B_1\setminus E_k\mbox{: }v_m(x)\neq 0\}\bigr|\leq 2W_k\Bigl({S\over m^p}\Bigr).
\end{equation}
In fact, from the $L^1$-norm lower semicontinuity of the map 
$$u\to\bigl|\{x\in B_1\setminus E_k\mbox{: } |u|(x)\neq 0\}\bigr|
\,=\,\int_{B_1\setminus E_k}\chi_{(0,\infty)}(|u|(x))dx,$$
it follows that
\begin{eqnarray*}
 m^p\bigl|\{x\in B_1\setminus E_k\mbox{: }v_m(x)\neq 0\}\bigr|
 & \leq & \liminf_{h\to\infty}m^p\bigl|
\{x\in B_1\setminus E_k\mbox{: }(v_{h,m}-u_h)(x)\neq 0\}\bigr|\\
\salt
& = &\liminf_{h\to\infty}\,m^p\bigl|E_{h,m}\setminus E_k\bigr| 
\,\leq \,2W_k\Bigl({S\over m^p}\Bigr).
\end{eqnarray*}
Now, setting $A_m:=\{x\in B_1\setminus E_k\mbox{: }v_m(x)\neq 0\}$, 
we  obtain from (\ref{seminf5}) that
\begin{equation}\label{seminf6}
\liminf_{h\to\infty}\int_{B_1} f_h(x,\E u_h)\dx
\geq \int_{B_1\setminus  (E_k\cup A_m)} f(0)\dx-\,CW_k\Bigl({S\over m^p}\Bigr)
\end{equation}
So, passing to the limit in (\ref{seminf5}) as 
$m\to\infty$ and using (\ref{estim4}) we finally obtain (\ref{seminf1})
 and this achieves the proof of the proposition.
\end{proof}
\vskip .1truecm
Now we are in the position to prove our main result.
\vskip .2truecm
\begin{proof}[The proof of Theorem \ref{lower1}.] Let $(u_h)$ be a sequence such that
$u_h$ converges strongly to $u$ in $L^1(\Om,\R^n)$ and $|E^ju_h|$ converges weak $*$ to 
the measure $\nu$ singular with respect to the Lebesgue measure. We assume that
$$\liminf_{h\to\infty}\int_\Om f(x,\E u_h)\dx
\,=\,\lim_{h\to\infty}\int_\Om f(x,\E u_h)\dx\,<\,\infty.$$
So, up to a subsequence, the sequence of  measures 
$f_h(x,\E u_h)\mathcal L^n\res\Om$ converges weakly $*$ to a positive 
measure $\mu$.
To prove (\ref{lower2}), it is enough to prove that
\begin{equation}\label{estim5}
\frac{d\mu}{d\mathcal L^n}(x_0)\,
\geq\,f(x_0,\mathcal Eu(x_0))\qquad\mbox{ a.e. }x_0\in\Om.
\end{equation}
In fact, from the lower semicontinuity of the total 
variations of measure with respect to weak $*$ 
convergence and from the inequality
(\ref{estim5}) it follows that
$$\liminf_{h\to\infty}\int_\Om f(x,\E u_h)\dx\geq 
\mu(\Om)\geq\int_\Om\frac{d\mu}{d\mathcal L^n}(x)\dx
\geq\int_\Om f(x,\mathcal Eu(x))\dx.$$
So, let us prove that (\ref{estim5}) holds. 
To this aim, we use a characterization of Carath\'eodory 
functions by Scorza-Dragoni (see e.g. \cite[Page 235]{ektem}), 
to get for every $i\in\N$ a compact set $K_i\subset\Om$ 
such that $|\Om\setminus K_i|<1/i$ and
$f|_{K_i\times\msim}$ is continuous in $K_i\times\msim$. 
Let $K_i^1$ be the set of Lebesgue points of the function $\chi_{K_i}$. We set
$$F:=\bigcup_{i\in\N}(K_i\cap K_i^1)$$
 and  it follows that $|\Om\setminus F|\,
\leq\,|\Om\setminus (K_i\cap K_i^1)|\,=\,
|\Om\setminus K_i|\,\leq\,1/i\to 0\quad\mbox{ as } i\to\infty$.\\
Let us fix $x_0\in F$ such that:
\begin{itemize}
\item[(i)] $x_0$ is an approximate differentiability 
point of $u$ and such that $\mathcal Eu(x_0)=\frac{\nabla u(x_0)+\nabla u(x_0)^T}{2}$;
\item[(ii)] $\displaystyle\frac{d\nu}{d\mathcal L^n}(x_0)
\,=\,\lim_{\e\to 0}\frac{\nu (B(x_0,\e))}{|B(x_0,\e)|}=0$;
\item[(iii)] $\displaystyle\frac{d\mu}{d\mathcal L^n}(x_0)
\,=\,\lim_{\e\to 0}\frac{\mu (B(x_0,\e))}{|B(x_0,\e)|}\,<\,\infty$.
\end{itemize}
Now, we consider a sequence $\e_k\searrow 0^+$ such that 
$\nu (\partial B(x_0,\e_k))=0$ and $\mu (\partial B(x_0,\e_k))=0$. 
Note that such a sequence
exists since $\{\e>0\mbox{: }\nu (\partial B(x_0,\e))>0\mbox{, }\nu (\partial B(x_0,\e))>0\}$ 
is at most a countable set. \\
From the approximate 
differentiability of $u$ at $x_0$ and the fact that
$u_h\to u$ strongly in $L^1(\Om,\R^n)$  we get 
\begin{equation}\label{estim6}
\lim_{k\to \infty}\lim_{h\to \infty}||u_{k,h}-w_0||_{L^1(\B_1,\R^n)}=0
\end{equation}
where 
$$u_{k,h}:=\frac{u_h(x_0+\e_ky)-u(x_0)}{\e_k}
\quad\mbox{ and }\quad w_0(y):=\nabla u(x_0)y.$$
 We have also that
\begin{eqnarray*}
|E^ju_{k,h}|(B_1)
&=&\int_{B_1\cap J_{u_{k,h}}}\bigl|(u^+_{k,h}-u^-_{k,h})\odot\nu_{u_{k,h}}\bigr|d\mathcal H^{n-1}\\
& = &\e_k^{-n}\int_{B(x_0,\e_k)\cap J_{u_h}}\bigl|(u^+_h-u^-_h)\odot\nu_{u_h}\bigr|d\mathcal H^{n-1}\\
& = &\frac{|E^ju_h|(B(x_0,\e_k))}{\e_k^n}\,\leq\, \frac{|E^ju_h|(\overline{B(x_0,\e_k)})}{\e_k^n}.
\end{eqnarray*}
Hence
\begin{eqnarray}\label{estim7}
\limsup_{k\to\infty}\limsup_{h\to\infty}|E^ju_{k,h}|(B_1) &\leq & 
\limsup_{k\to\infty}\limsup_{h\to\infty}\frac{|E^ju_h|(\overline{B(x_0,\e_k)})}{\e_k^n}\\
\nonumber &\leq & \limsup_{k\to\infty}\frac{\nu(\overline{B(x_0,\e_k)})}{\e_k^n}\,=\,0.
\end{eqnarray}
On the other hand, setting $f_k(y,\xi):=f(x_0+\e_ky,\xi)$ we get that
\begin{eqnarray*}
\frac{d\mu}{d\mathcal L^n}(x_0) 
&\geq &\limsup_{k\to\infty}\frac{\mu(\overline{B(x_0,\e_k)})}{|B(x_0,\e_k)|}\\
&\geq &\limsup_{k\to\infty}\limsup_{h\to\infty}\frac{1}{|B(x_0,\e_k)|}\int_{B(x_0,\e_k)}f(x,\E u_h)\dx\\
&\geq &\limsup_{k\to\infty}\limsup_{h\to\infty}\frac{1}{w_n}\int_{B_1}f(x_0+\e_ky,\E u_{k,h})\,dy\\
&= &\limsup_{k\to\infty}\limsup_{h\to\infty}\frac{1}{w_n}\int_{B_1}f_k(y,\E u_{k,h})\,dy.
\end{eqnarray*}
Therefore, by a standard diagonalization argument we may 
extract a subsequence $v_k:=u_{k,h_k}$ such that
$$
\lim_{k\to \infty}||v_k-w_0||_{L^1(B_1,\R^n)}=0\mbox{, }
\quad \lim_{k\to \infty}|E^jv_k|(B_1)=0$$
and
$$
\frac{d\mu}{d\mathcal L^n}(x_0)\,\geq\,
\limsup_{k\to\infty}\frac{1}{w_n}\int_{B_1}f_k(y,\E v_k)\,dy.
$$
 Now, since $x_0\in F$, there exist $i_0\in\N$ 
such that $x_0\in K_{i_0}\cap K^1_{i_0}$. So,
the sequence $\chi_{\frac{K_{i_0}-x_0}{\e_k}}$ 
converges strongly to $1$ in $L^1(B_1)$ and hence, up to a subsequence
$\chi_{\frac{K_{i_0}-x_0}{\e_k}}(y)\to 1$ for a.e. $y\in B_1$.
So, for $k$ large enough we have that  
$x_0+\e_ky\in K_{i_0}$ for a.e. $y\in B_1$. Hence, for every $\xi\in\msim$
we get  $\lim_{k\to\infty}f(x_0+\e_ky,\xi)\,=
\,f(x_0,\xi)$ for a.e. $y\in B_1$. 
Therefore, we  get
for a.e. $y\in B_1$,
\begin{equation}\label{estim8} 
\lim_{k\to\infty}f_k(y,\xi)\,=\,f(x_0,\xi)
\end{equation}
locally uniformly in $\msim$. So, applying Proposition \ref{step1} 
to the sequence $(v_k)$, we get that
$$\frac{d\mu}{d\mathcal L^n}(x_0)\,\geq\,
\liminf_{k\to\infty}\frac{1}{w_n}\int_{B_1}f_k(y,\E v_k)\,dy\,
\geq\,\frac{1}{w_n}\int_{B_1}f(x_0,\E u(x_0))\,dy\,=\,f(x_0,\E u(x_0))
$$
which gives (\ref{estim5}) and achieves the proof of the theorem.
\end{proof}
\section{Some examples and remarks}
 In the proof of Theorem \ref{lower1}, the assumption on 
$|E^ju_h|$ has played a crucial role in order to perform 
the blow-up argument. Note that any sequence $(u_h)\subset W^{1,p}(\Om,\R^n)$ such that 
$u_h\to u$ strongly in $L^1(\Om,\R^n)$ satisfies trivially the assumptions
of the theorem. For examples of sequences which are not necessarly
in $ W^{1,p}(\Om,\R^n)$, we consider here a variational problem 
 with  a uniform $L^\infty$ constraint on the admissible functions and 
 a unilateral constraint on their jump sets.\\ 
Let us  recall here the compactness criterion 
in $SBD$ by Bellettini-Coscia-Dal Maso \cite{BCDM}.

\begin{theorem}\label{Compactness-SBD}
 Let $\phi :[0,+\infty [\to [0,+\infty [$ be a non-decreasing 
function such that
\begin{equation}\label{Compactness-SBD1}
\lim\limits _{t\to +\infty }{{\phi (t)}\over t}=+\infty.
\end{equation}
Let $(u_h)$ be a sequence in $SBD(\Om )$ such that
\begin{equation}\label{Compactness-SBD2}
\int _\Om\vert u_h\vert dx +\vert
E^ju_h\vert (\Om )+\int _\Om\phi (\vert {\mathcal E}u_h\vert )dx+{\mathcal
  H}^{n-1}(J_{u_h})\leq C
\end{equation}
for some positive constant $C$ independent of $h$.\\
Then there exists a subsequence, still denoted by $(u_h)$ and a
function $u\in\sbd$ such that 
\begin{equation}\label{Compactness-SBD3}
u_h\to u \mbox{ strongly in }L^1_{\rm loc}(\Om ,\R^n),
\end{equation}
\begin{equation}\label{Compactness-SBD4}
{\mathcal E}u_h\rightharpoonup {\mathcal E}u\mbox{ weakly }\mbox{
  in }L^1(\Om ,\msim),
\end{equation}
\begin{equation}\label{Compactness-SBD5}     
E^ju_h\rightharpoonup E^ju\mbox{ weakly }\star\mbox{
  in }{\mathcal M}_b(\Om ,\msim),
\end{equation}
\begin{equation}\label{Compactness-SBD6}     
{\mathcal H}^{n-1}(J_u)\leq\liminf\limits _{h\to +\infty }{\mathcal
  H}^{n-1}(J_{u_h}).
\end{equation} 
\end{theorem}
In the next example we consider a variational problem for which the 
minimizing sequences satisfy the assumption on $|E^ju_h|$ in Theorem 
\ref{lower1}.
\begin{example}\label{exp0}
{\rm Let $K\neq\emptyset$ be a non closed subset of $\Om$ such that 
$0<\mathcal H^{n-1}(K)<\infty$ and  let $\{F(x)\}_{x\in\Om}$ be a family of 
uniformly bounded closed subsets of $\R^n$.
We consider the following variational 
problem:
\begin{equation}\label{min-const1}
\min_{\substack{u\in\sbd\\
 {\substack{J_u\subset K\\ u(x)\in F(x)\mbox{ a.e. in }\Om}}}}
\hskip -.6truecm \int_\Om f(x,\mathcal Eu)\dx
\end{equation}
with the function $f$ satisfying the assumptions of Theorem \ref{lower1}. 
Let us prove that Problem (\ref{min-const1}) admits a solution. \\
By the rectifiability of jump sets
 of BD functions, the inclusion $J_u\subset K$ will be 
intended up to a $\mathcal H^{n-1}$-negligible set.

Let  $(u_h)\subset
 \sbd$ be a minimizing sequence for problem (\ref{min-const1}). 
 By the assumptions, there exists $M>0$ such that  $||u_h||_\infty\leq M$
\begin{equation}\label{jumpmeas}
|E^ju_h|(\Om)\leq 2||u_h||_\infty\mathcal H^{n-1}(J_{u_h})
\leq 2M\mathcal H^{n-1}(K)<\infty
\end{equation}
and hence by the growth assumptions of $f$, (\ref{Compactness-SBD2}) is 
satisfied with $\phi(t)=t^p$. By Theorem \ref{Compactness-SBD},
 $u_h$ converges (up to a subsequence)
 strongly in $L^1(\Om,\R^n)$ to some function $u\in\sbd$.
Hence we get also $u(x)\in F(x)$ a.e. $x\in\Om$.

On the other hand, the sequence
 $|E^ju_h|$ converges (up to a subsequence) weakly $\star$
 to some positive measure $\nu$. It easily follows from (\ref{jumpmeas}) 
that the measure $\nu$
 is concentrated on the set $K$. Therefore $\nu$ is singular with respect 
to the Lebesgue measure. So, by Theorem \ref{lower1}, we have that
$$\int_\Om f(x,\mathcal Eu)\dx\leq\liminf_{h\to\infty}
\int_\Om f(x,\mathcal Eu_h)\dx.$$
Now let us prove that $u$ verifies the constraint $J_u\subset K$ up to a 
$\mathcal H^{n-1}$-negligible set. This is obtained by slicing.
To this aim we recall the notations 
for one-dimensional sections of BD functions. 

Given $\xi\in\R^n$ with $\xi\neq 0$, we set
$$\pi^\xi:=\{y\in\R^n\mbox{: }(y,\xi)=0\}$$
and for every $y\in\pi^\xi$ and for every $B\in\mathcal B(\Om)$,
$$B_y^\xi:=\{t\in\R\mbox{: }y+t\xi\in B\}\quad\mbox{ and }
\quad B^\xi:=\{y\in\pi^\xi\mbox{: }B^\xi_y\neq\emptyset\}.$$
For every $u\in L^1(\Om,\R^n)$ we set 
$$u_y^\xi(t):=(u(y+t\xi),\xi).$$
It has been proved in \cite{ACDM} that, if $u\in\sbd$ then for 
$\mathcal H^{n-1}$-a.e. $y\in\Om^\xi$, $u_y^\xi\in\,\,$SBV$(\Om_y^\xi)$.
Viceversa, assume that 
$$u_y^\xi\in\mbox{SBV}(\Om_y^\xi)\mbox{ for } 
\mathcal H^{n-1}\mbox{-a.e. }y\in\Om^\xi\quad\mbox{and}\quad
 \int_{\Om^\xi}|Du^\xi_y|(\Om^\xi_y)d\mathcal H^{n-1}(y)<\infty$$ 
for every $\xi=\xi_i+\xi_j$, $i,j=1,\cdots,n$, with $(\xi_i)_{i=1}^n$ being 
an orthonormal basis in $\R^n$. Then $u\in\sbd$.\\
Setting $J^\xi_u:=\{x\in J_u\mbox{: }(u^+(x)-u^-(x),\xi)\neq 0\}$, it follows
 from Fubini's theorem that 
\begin{equation}\label{jumpxi}
\mathcal H^{n-1}\bigl(J_u\setminus J_u^\xi\bigr)=0\quad\mbox{ for }
\mathcal H^{n-1}\mbox{-a.e. } 
\xi\in\mathcal S^{n-1}.\end{equation}
From the structure theorem for BD 
functions (see \cite[theorem 5.1]{ACDM}) we have also  
$$J_{u_y^\xi}=\bigl(J_u^\xi\bigr)_y^\xi\quad\mbox{ for a.e. }
 y\in\Om^\xi.$$

Now we can prove that the limit $u$ of the minimizing sequence $(u_h)$ for 
Problem (\ref{min-const1}) satifies the constraint $J_u\subset K$. \\
Let
 $\xi\in\mathcal S^{n-1}$ be such that (\ref{jumpxi}) holds.
Following the proof of Theorem \ref{Compactness-SBD}, we get that 
the sequence of one-dimensional section $(u_{h,y}^\xi)$ of the minimizing 
 sequence $(u_h)$ satisfies  
 the assumptions of the SBV compactness theorem and 
from $J_{u_h}\subset K$ we get 
$$J_{u_{h,y}^\xi}=(J_u^\xi)_y^\xi\subset K_y^\xi\mbox{ with }
\mathcal H^0(K_y^\xi)<\infty\mbox{ for }\mathcal H^{n-1}\mbox{-a.e. }y\in\Om^\xi.$$
Therefore the limit function $u_y^\xi$ has also its
 jump set contained in the finite set  $K_y^\xi$. In fact, it is easy to see
that the jump set $J_{u_y^\xi}$ is contained in the set of limits of the jump
 points of $u_{h,y}^\xi$. 
 Now from 
$(J_u^\xi)_y^\xi=J_{u_y^\xi}\subset K_y^\xi
\mbox{ for }\mathcal H^{n-1}\mbox{-a.e. }y\in\Om^\xi$, we get 
$J_u^\xi\subset K$ up to a $\mathcal H^{n-1}$-negligible set and hence by
 (\ref{jumpxi}), also
$J_u\subset K$ up to a $\mathcal H^{n-1}$-negligible set.
\qed
} 
\end{example}
\begin{remark}\label{rmK}
{\rm Note that the set $K$ has been taken non closed in order to avoid
 the easy case where the minimizing sequences $(u_h)$ and their limit $u$ belong
to the space 
$$LD(\Om\setminus K)
:=\{u\in L^1(\Om\setminus K,\R^n)\mbox{: }Eu\in L^1(\Om\setminus K,\msim)\}$$
for which the lower semicontinuity of the functional
$$\int_\Om f(x,\mathcal Eu)\dx=\int_{\Om\setminus K} f(x,\mathcal Eu)\dx$$
in the strong topology of $L^1(\Om\setminus K,\R^n)$ follows from \cite[Theorem 3.1]{EBO2}. 
}
\end{remark}
\vskip .2truecm
As we have seen in the previous example, the minimizing sequences for problem (\ref{min-const1})
  satisfy the assumptions of both Theorems \ref{lower1} 
and \ref{Compactness-SBD}. 
However, unlike the assumptions (\ref{jumpset}) 
in Theorem \ref{theoamb} and (\ref{kristjump}) in \cite{KRIS}, 
which are consistent with the compactness criterion in $SBV$ 
 (see for instance \cite[Theorem 4.8]{AMFUPA}), the assumption of Theorem \ref{lower1} on the measures 
 $|E^ju_h|$ is not always compatible with the compactness criterion in Theorem 
\ref{Compactness-SBD}.

In the following example, we construct a sequence $(u_h)\subset SBD(\Om)$ 
which satisfies the compactness criterion in $SBD$
  while $|E^ju_h|$ converges to a  measure proportional 
 to the Lebesgue measure. 
\begin{example}\label{examp1}
{\rm We consider in $\R^2$ the open squares
$$\Om:=\bigl(0,2\bigr)\times\bigl(0,2\bigr)\quad\mbox{ and }\quad
\Om_h:=\Bigl(\frac{1}{h}-\frac{1}{h^2},\frac{1}{h}+\frac{1}{h^2}\Bigr)
\times\Bigl(\frac{1}{h}-\frac{1}{h^2},\frac{1}{h}+\frac{1}{h^2}\Bigr).$$
We set
$$
E_h:=\bigcup_{(i,j)\in I_h\times I_h}\bigl(\Om_h +(i,j)\bigr)\quad\mbox{ with }
\quad I_h:=\{0,2/h,4/h,\cdots,2-2/h\}.$$
 Let $(u_h)$ be the sequence defined by $u_h:=(\chi_{E_h},0)$ and let $E_{h,i,j}:=\Om_h +(i,j)$. 
 By easy computations we get 
$$E^ju_h=Eu_h=\sum_{(i,j)\in I_h\times I_h}(1,0)\odot\nu_{E_{h,i,j}}\mathcal H^1\res\partial E_{h,i,j}$$
where $\nu_{E_{h,i,j}}$ 
is the unit normal vector to $\partial E_{h,i,j}$.
 
We have that $|E^ju_h|(\Om)=2\sqrt 2+4$ and \mbox{$\mathcal H^1(J_{u_h}\cap\Om)=8$}. Thus,
 the sequence $(u_h)$
 satisfies the assumptions of Theorem \ref{Compactness-SBD}. \\
However, the sequence $|E^ju_h|$ converges weakly $*$ to the measure 
\mbox{$(\sqrt 2+2)\mathcal L^2\res\Om$.} Indeed, let $M^{i,j}_h,\,N^{i,j}_h$ and $L^{i,j}_h,\,K^{i,j}_h$ be
respectively  the two vertical  and horizontal sides of 
the square $E_{h,i,j}$.
 It is easy to see that
\begin{equation}\label{J1}
|E^ju_h|=\sum_{(i,j)\in I_h\times I_h}\Bigl(\mathcal H^1\res M^{i,j}_h+\mathcal H^1\res N^{i,j}_h
+\frac{\sqrt 2}{2}\mathcal H^1\res L^{i,j}_h+\frac{\sqrt 2}{2}\mathcal H^1\res K^{i,j}_h\Bigr).
\end{equation}
Now let $\varphi\in C_c(\Om)$. It is easy to see that
\begin{equation}\label{J2}
\lim_{h\to\infty}\sum_{(i,j)\in I_h\times I_h}\int_{S^{i,j}_h}
\varphi d\mathcal H^1= \int_\Om\varphi\dx\quad\mbox{ for }
S^{i,j}_h=M^{i,j}_h,N^{i,j}_h,K^{i,j}_h,L^{i,j}_h.
\end{equation}
Therefore from (\ref{J1}) and (\ref{J2}) we get
$$\lim_{h\to\infty}\int_\Om\varphi \,d|E^ju_h|=(\sqrt 2+2)\int_\Om\varphi\dx.$$
\qed}
\end{example}
\vskip .4truecm
\centerline {\sc Acknowledgements:}
 \vskip .1truecm
\noindent The author is grateful to Luigi Ambrosio for helpful 
discussions on the subject of the paper. This work is supported by the
National Research Foundation of South Africa.

\end{document}